\documentclass[12pt]{article}
\usepackage{graphicx}
\usepackage{amssymb}
\usepackage{amsfonts}
\usepackage{amsthm}
\usepackage{amsfonts}

\newcommand \ra {\rightarrow}

\newcommand{\ba}[1]{\begin{array}{#1}}
\newcommand{\ea}{\end{array}}

\newcommand{\be}{\begin{equation}}
\newcommand{\ee}{\end{equation}}
\newcommand{\bea}{\begin{eqnarray}}
\newcommand{\eea}{\end{eqnarray}}
\newcommand{\beann}{\begin{eqnarray*}}
\newcommand{\eeann}{\end{eqnarray*}}

\bibstyle{ams}

\begin{document}

\title{A non-intersecting random walk on the Manhattan lattice and SLE$_6$}

\author{Tom Kennedy
\\Department of Mathematics
\\University of Arizona
\\Tucson, AZ 85721
\\ email: tgk@math.arizona.edu
}

\maketitle 

\begin{abstract}
We consider a random walk on the Manhattan lattice. The walker must follow
the orientations of the bonds in this lattice, and the walker is not allowed
to visit a site more than once. When both possible steps are allowed, the 
walker chooses between them with equal probability. 
The walks generated by this model are known to be related to interfaces 
for bond percolation on a square lattice. So it is natural to 
conjecture that the scaling limit is SLE$_6$. We test this conjecture with 
Monte Carlo simulations of the random walk model and find strong support
for the conjecture. 
\end{abstract}

\section{Introduction}

There are many different models of random walks on a lattice 
which generate walks which do not intersect themselves or 
in which self-intersections are disfavored. 
In the model that is usually referred to as ``the self-avoiding walk'', one
considers all nearest neighbor walks of length $N$ that do not have any 
self-intersections and defines the probability measure to be the uniform
measure on this set of walks. One would then like to let 
$N \rightarrow \infty$. In this model there is not a simple relation between 
the walks with $N+1$ steps and those with $N$ steps. 
There are variations on this model. One can consider a bond avoiding random walk
in which the walk is allowed to self-intersect, 
but is not allowed to traverse any bond more than once. Or one can 
consider a weakly self-avoiding walk in which all nearest neighbors 
random walks with $N$ steps are allowed, and the probability of a walk 
is proportional to $e^{- \beta I}$ where $I$ is the number of self-intersections
and $\beta>0$ is a parameter. 
See \cite{madras_slade} for more on these models. 
It is conjectured that all these models 
have the same scaling limit and that the limiting process is SLE$_{8/3}$ 
\cite{lsw_saw}. Simulations of the self-avoiding walk support this 
conjecture \cite{Kennedya,Kennedyb}. 
Another model which generates walks with no self intersections is the 
loop-erased random walk (LERW). One takes an ordinary nearest neighbor
random walk on the lattice and erases the loops it forms in chronological
order. This model has been proved to converge to SLE$_2$ in the scaling 
limit \cite{lerw_sle}.

In the models above, there is no algorithm that will grow a 
random walk with the given probability distribution. There is 
another class of random walks without self-intersections 
in which the probability measure
is defined dynamically, i.e., there is an algorithm that generates a sample 
of the walk one step at a time.
The simplest model is to let the walk choose its next step by 
randomly picking one of its unoccupied nearest neighbors with equal 
probability. The problem with this model is that there may not be any 
unoccupied nearest neighbors; the walk can get trapped. This trapping can 
be avoided by letting the walk choose any nearest neighbor with 
probabilities that favor those sites that have not been visited yet. 
Models of this type are typically called the true self-avoiding walk 
or the myopic self-avoiding walk \cite{amit_et_al}. 
Another model is to let the walk 
pick with equal probability one of the nearest neighbors that satisfies two 
conditions - the walk has not visited the neighbor before and there is an 
infinite path from the neighbor that avoids sites that have been 
visited before. 
This model was introduced in the physics literature in the 80's under two 
names - the smart kinetic walk and the infinitely growing self-avoiding walk
\cite{kremer1985IGSAW,weinrib_trugman}.
On the hexagonal lattice it is equivalent to 
the percolation explorer and so is known to have SLE$_6$ as its scaling limit. 
\cite{smirnov,camia_newman}.
There is numerical evidence that the scaling limit for the model on 
the square and triangular lattices is also SLE$_6$ \cite{dai,Kennedyc}.

In this paper we study a random walk model on the Manhattan lattice 
in which the walk is not allowed to visit a site more than once. 
The Manhattan lattice is an oriented square lattice in which the orientations
are constant along horizontal and vertical lines and alternate as we 
move up or down from a horizontal line or right or left from a vertical line. 
So if the walk is at $\omega(n)$ at time $n$, there are at most 
two possibilities for $\omega(n+1)$. If both of the possibilities have not 
been visited before, we randomly choose one with equal probability. 
If just one possibility has not been visited before, $\omega(n+1)$ is taken to 
be that site. If both of the possibilities have been visited before, 
then the walk will be trapped. It was observed long ago that this can only 
happen when $\omega(n)$ is a nearest neighbor of the starting point of the 
walk \cite{HH}. As we will discuss in the next section, this argument also 
shows that if we restrict the walk to certain domains, the 
walk will never be trapped. The particular domains we will use for our
simulations are a slit plane, the upper right quadrant and the complement of
the upper right quadrant. 

As we will discuss in the next section, this random walk model is related 
to interfaces for bond percolation on a square lattice \cite{gunn_ortuno}.
So a natural conjecture for the scaling limit of this model is SLE$_6$. 
We will test this conjecture with numerical simulations of the model.

\section{Definition and equivalent forms of the model}

Our model is a random walk on the Manhattan lattice which is not allowed
to visit a site more than once. When both possible directions for a step 
have not been visited before, the walk chooses one with equal probability.
If both of the two possible sites for the next step of the walk have been 
visited before, the walk will be trapped. Hemmer and Hemmer \cite{HH}
showed that this can only happen when the site is a nearest neighbor of the 
starting point of the walk. We will repeat their argument. 

\begin{figure}[tbh]
\includegraphics{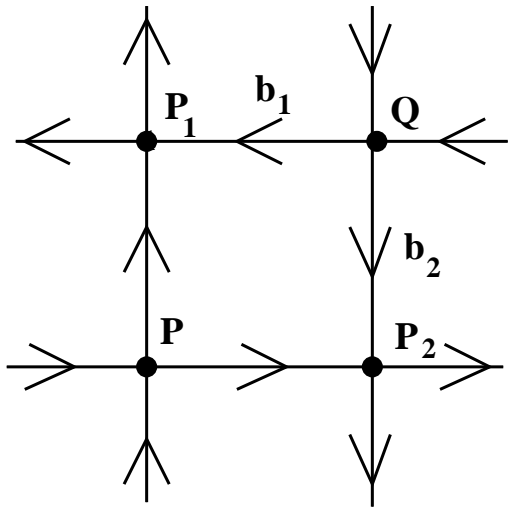}
\caption{The figure illustrates the argument that the walk cannot get 
trapped at $P$.
}
\label{non_trapping}
\end{figure}

In figure \ref{non_trapping}  site $P$ is the current location of the walk. 
We assume the walk does not start at $P_1$ or $P_2$. So if $P_1$ has 
been visited before, the walk must have followed bond $b_1$. 
And if $P_2$ has been visited before, 
the walk must have followed bond $b_2$. But this implies that 
site $Q$ was visited twice, which is a contradiction. 

\begin{figure}[tbh]
\includegraphics{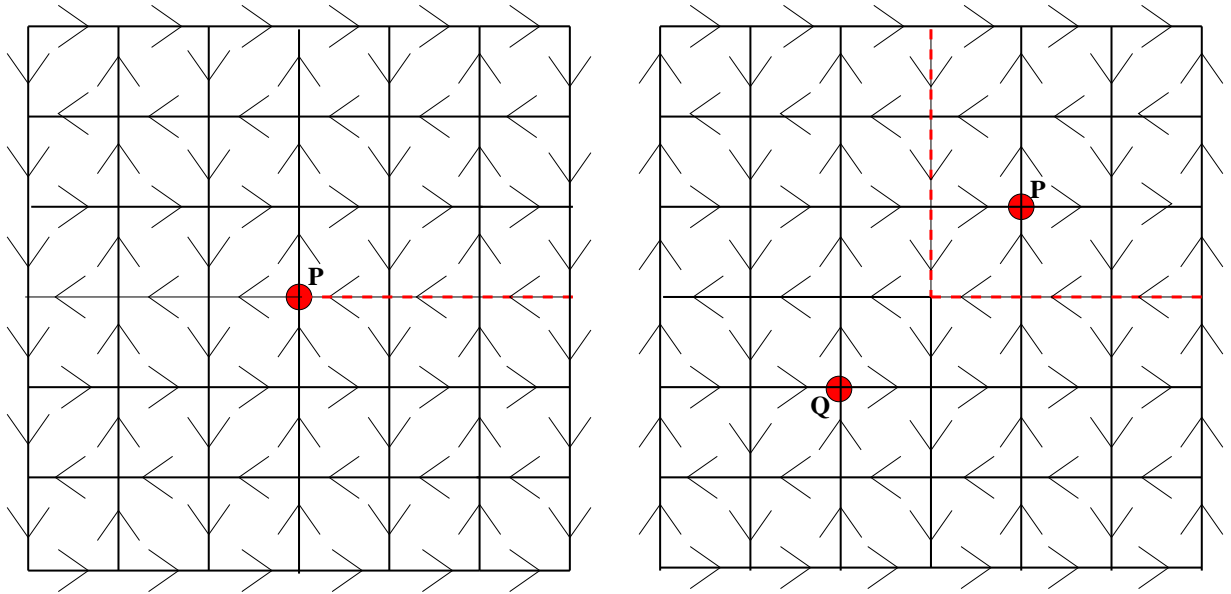}
\caption{The random walk on the Manhattan lattice with no self-intersections
is not allowed to visit sites on the dashed (red) lines. 
For these three domains, if it starts at $P$ or $Q$, it will never be trapped.
}
\label{domains}
\end{figure}

Using results from percolation, it can be proved that the walk will return 
to a nearest neighbor of the starting point and be trapped 
with probability one \cite{Grimmett}. 
By restricting the walk to certain domains, we can ensure that it will never
be trapped. Figure \ref{domains} gives three examples. 
In these figures the walk 
is not allowed to visit sites on the boundary which is given by the 
dashed (red) line(s). In the left domain in figure \ref{domains}, 
the walk is started at $P$. 
In the right domain, we can start the walk at $P$ in which case the 
walk will remain in the upper right quadrant, 
or we can start it at $Q$ in which case 
it will remain in the complement of the upper right quadrant. 
We will refer to these three domains as the slit plane, 
the 90 degree wedge and the 270 degree wedge. 
It is easy to check using arguments similar to the previous paragraph 
that the walk will never be trapped in these domains with these starting 
points. 

\begin{figure}[tbh]
\includegraphics{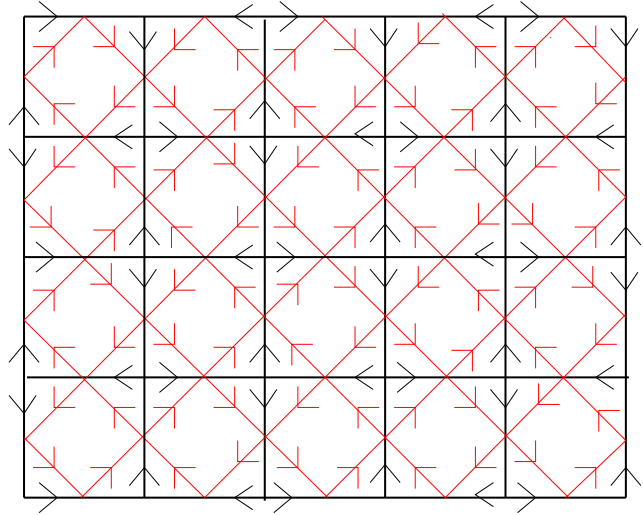}
\caption{The L lattice is the lattice whose bonds are horizontal and vertical. 
The associated Manhattan lattice is the (red) lattice whose bonds
are at 45 degrees with respect to horizontal and vertical.
}
\label{covering}
\end{figure}

There is an equivalent formulation of the random walk model we have 
been considering. 
It uses a different oriented square lattice which is 
sometimes called the L lattice. The Manhattan lattice is the ``covering 
lattice'' for the L lattice in the sense of Kasteleyn \cite{kasteleyn}. 
Figure \ref{covering} shows the L lattice and the associated Manhattan lattice.
The random walk model on the L lattice is defined as follows. 
It follows the orientations in the lattice, so at each step it can 
only turn right or left. 
It is not allowed to traverse a particular bond more than once, but it 
can visit a site more than once (but only twice given the bond constraint). 
When both turns are possible, it chooses between them with equal probability.
Given such a random walk on the L lattice, the midpoints of the bonds
in the walk will be a walk on the Manhattan lattice that never 
visits a site in the Manhattan lattice more than once. 
This gives a 1-1 correspondence between $n$-step walks on the 
L lattice with no repeated bonds and $(n-1)$-step random walks 
on the Manhattan lattice with no self-intersections. 
The equivalence of these two random walk models on the Manhattan lattice 
and the L lattice was observed by Malakis \cite{malakis}.

We note that the term ``self-avoiding walk on the Manhattan lattice'' 
is used to refer to the model defined by putting the uniform probability
measure on the set of self-avoiding walks on the Manhattan lattice with a 
fixed number of steps. (For example, \cite{malakis} considers this model.) 
This is a completely different model from the 
model we consider in this paper. We have eschewed the adjective
``self-avoiding'' for our model to avoid confusion with this other model
on the Manhattan lattice.

This random walk on the L lattice that does not repeat bonds can also 
be formulated as a deterministic random walk in a random environment. 
This is a special case of the model considered by Gunn and Ortuno  
\cite{gunn_ortuno}. When the random walk makes a turn at a site it has
not visited before, we place a mirror at that site oriented so that the 
turn corresponds to the walk being reflected by the mirror. Note that 
when the walk returns to a site for a second time, the direction of the 
turn it must make is consistent with the orientation of the mirror that
is already at the site. So the mirrors have two possible orientations and 
they have equal probability. Rather than introduce the mirrors as the 
walk evolves, we can first randomly place mirrors on all the lattice sites.
The walk then evolves in a deterministic manner in this random 
environment. Note that we must generate a new random environment each time 
we want to generate a new sample of the random walk. 

Finally, we review the relation of the random walk on the L lattice
to interfaces in bond percolation on a square lattice.
It is worth noting at the start that the bond percolation process 
does not take place on the L lattice, but rather on a lattice that 
we will define shortly.
The sites in the dual lattice for the L lattice are at the centers of 
the squares in the L lattice. 
We take the length of the mirrors to be $\sqrt{2}$. 
(The lattice spacing is $1$.) Then the ends of the mirrors are at sites
in the dual lattice which are next nearest neighbors.
The dual lattice is a bipartite lattice, i.e., we can label 
the sites in the dual lattice as even or odd in such a way that an odd
site has only even sites as its nearest neighbors and vice versa.  
So the ends of the mirrors will either be both odd or both even.  Hence 
we can label the mirrors as even or odd. If we just consider the even sites in 
the dual lattice, we have a square lattice with spacing $\sqrt{2}$, rotated
by $45$ degrees with respect to the dual lattice. We will refer to 
this lattice as the even half-dual lattice. The odd half-dual lattice is 
defined similarly. The even mirrors are bonds in the even half-dual 
lattice, and the odd mirrors are bonds in the odd half-dual lattice.
Note that each site in the L lattice is the midpoint of one bond
in the even half-dual lattice and one bond in the odd half-dual lattice. 
There is a mirror at every site in the L lattice, and it will be the odd bond 
with probabilty $1/2$ or the even bond with probability $1/2$. 
So the even mirrors have the 
distribution of bond percolation on the even half-dual lattice. 
Likewise, the odd mirrors have the distribution of bond percolation on 
the odd half-dual lattice. These two percolations processes are not independent.
If fact, the configuration of odd mirrors completely determines the 
configuation of the even mirrors and vice versa. 

\begin{figure}[tbh]
\includegraphics{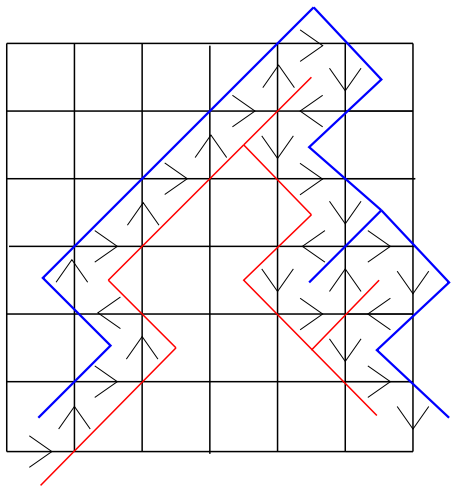}
\caption{A portion of a random walk on the L lattice is shown with arrows. 
It is sandwiched between a connected component of even mirrors
and a connected component of odd mirrors.}
\label{percolation}
\end{figure}

Gunn and Ortuno \cite{gunn_ortuno} showed that 
the random walk on the L lattice is sandwiched between a connected cluster 
for the bond percolation on the odd half-dual lattice and a connected 
cluster for the bond percolation on the even half-dual lattice. 
Figure \ref{percolation} shows a section of a random walk. 
Only the mirrors that touch the random 
walk are shown. As one traverses the random walk, one of the connected 
clusters is always on the right and the other is always on the left. 
So we can think of the random walk as either tracing out an interface
for the bond percolation process on the even half-dual lattice or an 
interface for the bond percolation process on the odd half-dual lattice.
Since our random walk on the Manhattan lattice with no self-intersections is 
equivalent to the random walk on the L lattice with no repeated bonds,
the random walk on the Manhattan lattice is also related to interfaces 
in bond percolation on the square lattice. 
This relation was described in \cite{bradley}.  
Ziff, Cummings and Stells studied other random walk models
that are related to percolation interfaces \cite{ziff1984}.

\section{Tests for SLE$_6$}

\subsection{Hitting distribution}

We test the conjecture that the scaling limit is SLE$_6$ by studying 
two quantities that can be computed explicitly for SLE$_6$. 
The first quantity is a type of hitting distribution. 
We first consider our model in the left domain in figure \ref{domains}, i.e., 
the plane slit along the positive real axis.  
We consider a circle of radius $r$ and stop the random walk 
when it first hits this circle. 
The random variable we study is the polar angle of 
the point where the walk first hits the circle. The distribution
of this random variable for SLE$_6$ is known explicitly. 

Our walk is going from $0$ to $\infty$, so in the scaling limit 
our conjecture is that we get SLE$_6$ in the slit
plane going from $0$ to $\infty$. We will refer to the subset of the slit
plane with $|z|<r$ as the slit disc.
Because of the locality property of 
SLE$_6$, up until the time the SLE$_6$ curve hits the circle of 
radius $r$, the distribution of the curve is the same as that of 
an SLE$_6$ curve in the slit disc going from $0$ to any fixed point 
on the circle. 
Note that each point on the positive real axis is really two boundary points
for the slit plane, one as we approach the axis from above and one as 
we approach it from below. We distinguish the two boundary points 
that correspond to an $x>0$ by $x^+$ and $x^-$.  
We take the terminal point for the SLE$_6$ to be $r^-$. 
So we consider SLE$_6$ in the 
slit disc with radius $r$ starting at $0$ and ending at $r^-$ and 
want to know the distribution of the first point where it hits the 
circle of radius $r$. We take the conformal map of the slit disc 
to the half plane which sends $0$ to $0$, $r^-$ to $\infty$ and $r^+$ to 
$1$. The boundary of the slit disc is mapped to the real axis, and 
the part of the boundary of the slit disc that is the circle of radius 
$r$ is mapped to $[1,\infty)$.
So we can find the distribution of the hitting point for the circle of radius 
$r$ if we can find the distribution of where SLE$_6$ in the half plane first
hits $[1,\infty)$. 

This distribution is known.
For $\kappa>4$, the SLE curve will touch the real 
axis infinitely often. Let $\gamma(t)$ denote the SLE trace, and let 
$t^*$ be the first time it touches the subset $[1,\infty)$ 
of the real axis.
So $\gamma(t^*)$ is the place it first hits $[1,\infty)$. 
Its distribution is the following; see proposition 6.34 
in \cite{lawler_book}. 
\bea
P(\gamma(t_*) < 1+x) = c I({x \over x+1}), \\
I(x)= \int_0^{x} u^{-2/3} (1-u)^{-2/3} \, du
\eea
with $c = 1/I(1)$. 

Returning to our walk in the slit disc, the conformal map of the slit disc 
to the half plane which sends $0$ to $0$, $r^-$ to $\infty$ and $r^+$ to
$1$ can be constructed as follows. We take $r=1$.
Define
\bea
\phi(z)= {1 \over 2} (z + {1 \over z}), \quad 
\psi(z)= {2 \over 1+z}
\eea
First we apply the map $z \ra \sqrt{z}$ to map the slit disc to the 
half disc in the upper half plane. 
Then we apply the map $\phi(z)$ to map this half disc to the half plane
below the real axis. Finally we apply the Moibius transformation $\psi(z)$ 
to get the points $0$, $r^+$ and $r^-$ to go to the appropriate points. 
So the overall conformal map is given by $\psi(\phi(\sqrt{z}))$.
We have
\bea
\psi(\phi(\sqrt{R e^{i \theta}})) = 1 + x, \quad
x=\frac{1-cos(\theta/2)}{1+cos(\theta/2)}
\eea
and so the cumulative distribution function (CDF) 
for the hitting distribution of the walk in the half plane is 
$F(\theta)=c I(x)$ with $x$ given above. 

For the 90 degree wedge, the conformal map $z \ra z^4$ maps the region to 
the slit plane, so the CDF is given by the above with the only change
being that we replace $\theta$ by $4 \theta$ in the definition of $x$. 
For our 270 degree 
wedge, the polar angle ranges from $\pi/2$ to $2 \pi$. 
The only change we need to make in the above is to replace 
$\theta$ in the definition of $x$ by $\frac{4}{3}(\theta-\pi/2)$.

Since we are looking at the first time that the walk hits a circle of radius
$R$, we need only simulate the walk up until the time it hits the circle. 
Unlike the next test that we will consider, we do not need to worry about 
the limit of taking the number of steps in the walk to infinity.
The simulations will be done with a lattice spacing of $1$, so 
$\delta=1/R$ can be thought of as the effective lattice spacing. 
The scaling limit is given by letting $\delta \ra 0$, i.e., $R \ra \infty$. 

Before we take the scaling limit, the random variable we are studying is 
discrete; there are only a finite number of points where the walk can first 
hit the circle. Of course, in the scaling limit it should converge to 
a continuous random variable. The discreteness of the random variable 
will be quite visible in our simulation results.
We can reduce the effect of this discreteness in the following way. 
As the radius $R$ changes, the set of possible values of the discrete random 
variable changes. So if we average the random variable over an interval of 
radii $R$, we get a continuous random variable.
We will refer to this as ``averaging over $R$.''
So if we let $\Theta_R$ denote that random variable for a circle of radius 
$R$, then we consider instead the averaged random variable 
\bea
\Theta_{R_0,R_1} = {1 \over R_1-R_0} \int_{R_0}^{R_1} \, \Theta_r \, dr
\eea
We will always take $R_0=R, R_1=2R$. 
So the scaling limit is given by taking $R \ra \infty$. 
Obviously, the distribution of the limit of $\Theta_{R,2R}$ should be the same
as the distribution of the limit of $\Theta_R$. 

\subsection{Pass right function}
     
The other quantity we study is the following. Let $z_0$ be a point in 
the domain where the walk is taking place. 
We consider the probability that the SLE trace passes 
to the right of the point. Since the SLE trace has self-intersection points
(as does the random walk), the definition of passing right of a point is not 
completely obvious. It can be defined using winding numbers as in 
\cite{schramm}. A more practical definition for computational purposes is
the following. Given the SLE  curve (or a random walk), take a curve 
from $z_0$ to a point on the positive real axis which is generic in the sense 
that it does not pass through any self-intersection points of the SLE curve 
(or the random walk). Count the number of intersections of this curve with 
the SLE curve (or the random walk). If this number is odd, the SLE curve 
(or random walk) passes right of $z_0$, if even it passes left. 

Schramm derived an explicit formula for this probability 
for SLE$_\kappa$  \cite{schramm}. In domains such as the half-plane, the 
slit-plane or the wedges we consider, by the dilation invariance of SLE 
this probability only depends on the polar angle of the point. 
We will refer to this probability of passing right of a point 
with polar angle $\theta$ as the pass-right function and denote it by 
$p_{SLE}(\theta)$. 
In the half plane, Schramm's formula is 
\bea
p_{SLE}(\theta) = c \int_0^\theta [\sin(t)]^{-2+8/\kappa} dt
\eea
where the constant $c$ is defined by $p(\pi)=1$. 
The half plane is mapped onto our three domains by a map of the form 
$z \rightarrow z^p$, so the formula for our domains is given by replacing 
$\theta$ in the right side by $\theta/2$ for the slit plane, by $2 \theta$
for the $90$ degree wedge and by $\frac{2}{3}(\theta-\pi/2)$ for the $270$ 
degree wedge. 

In our simulations we will compute the probability of passing right 
of $R e^{i\theta}$ with $R$ fixed and $\theta$ varying. We denote this
probability by $p_R(\theta)$. 
If the scaling limit of our random walk is SLE$_6$, 
$p_R(\theta)$ should converge to $p_{SLE}(\theta)$ as $R \rightarrow \infty$.  

Note that there is not a natural time at which we should stop the generation 
of the walk when we study the probability of passing right. 
No matter how far the walk is outside the circle of radius $R$, 
there is always some probability that the walk will cross this circle again
and so change whether some points on the circle are right or left of the walk. 
Unlike the first test we considered, we need to let the number of steps in 
the walk go to $\infty$. 
The simulation uses a lattice of spacing 1. 
If we rescale so that the circle has radius $1$, then 
the lattice spacing becomes $1/R$. So we will refer to $\delta=1/R$ as the 
lattice spacing. Note that how ``close'' $N$ is to infinity depends on $R$. 
The number of steps it takes the walk to first reach the circle of radius 
$R$ is of order $R^{1/\nu}$. (We take $\nu=4/7$.) So we define $n=N/R^{1/\nu}$. 
The proper way to measure how large $N$ is, is to consider how large $n$ is
with respect to $1$.  
We now must take a double limit to obtain the scaling limit; $\delta$ must
go to zero and $n$ must go to $\infty$.  

For finite $R$, the function $p_R(\theta)$ is a step function since there
are only a finite number of points where the walk can cross the circle 
of radius. As with the previous random variable, we can smooth out this 
function by averaging $R$ over some interval. So we define 
\bea
p_{R_0,R_1}(\theta) = {1 \over R_1-R_0} \int_{R_0}^{R_1} \, p_r(\theta) \, dr
\eea
We will always take $R_0=R, R_1=2R$. 
So if we define $n=N/R^{1/\nu}$, then the scaling limit is given by 
letting $R$ go to $\infty$ and $n$ go to $\infty$. 

\section{Simulation results}

\subsection{Hitting distribution}

\begin{figure}[tbh]
\includegraphics{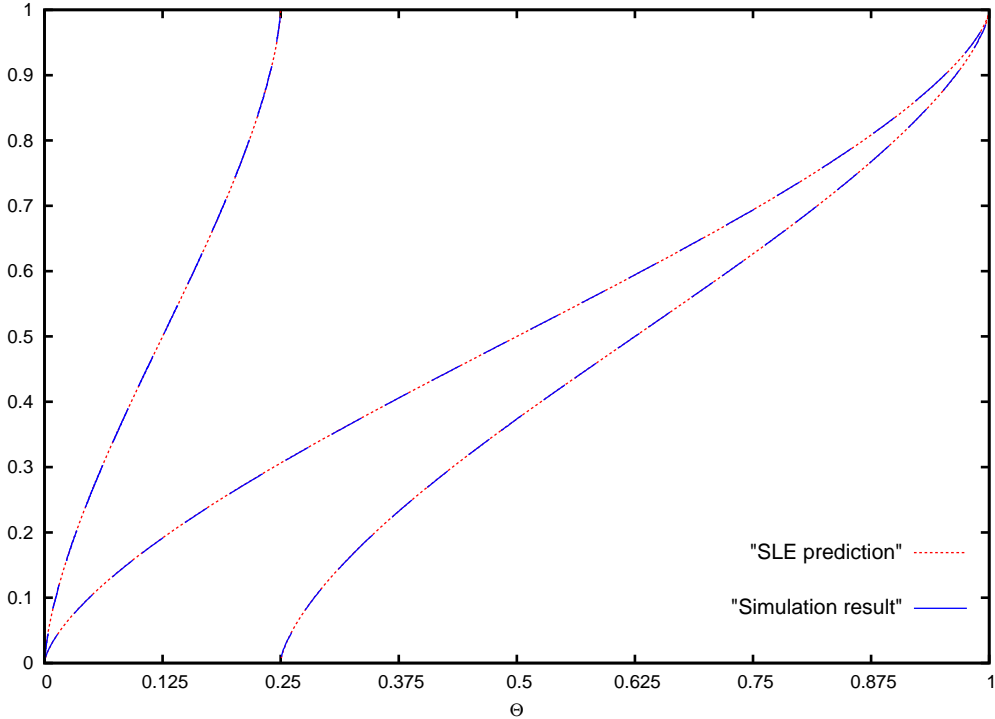}
\caption{Comparison of the CDF for the hitting distribution for 
the random walk and SLE$_6$ for three domains. 
The two curves for the slit plane have 
$\Theta \in [0,1]$, for the 90 degree wedge $\Theta \in [0,1/4]$ and for 
the 270 degree wedge $\Theta \in [1/4,1]$.}
\label{bad2_exit_cdf}
\end{figure}

We first consider the hitting distribution for three domains - 
the slit plane, the 90 degree wedge and the 270 degree wedge. 
Recall that by the hitting distribution we mean the polar angle of the 
point where the walk first hits a circle of radius $R$. 
The scaling limit is taken by letting $R \rightarrow \infty$. 
In all our figures the variable $\Theta$ is actually the polar angle 
divided by $2 \pi$. So it runs from $0$ to $1$ for the slit plane, from 
$0$ to $1/4$ for the $90$ degree wedge and from $1/4$ to $1$ for the 
$270$ degree wedge. 
In figure \ref{bad2_exit_cdf} we plot the CDF for this hitting distribution
for the three domains. For each domain two CDF's are plotted - 
the simulation CDF and the CDF predicted by SLE$_6$. 
They agree so well that the curves cannot be distinguished. 
We plot only ``half'' of each curve so that the overlapping curves can 
be seen. The simulations shown used a radius of $R=800$.

\begin{figure}[tbh]
\includegraphics{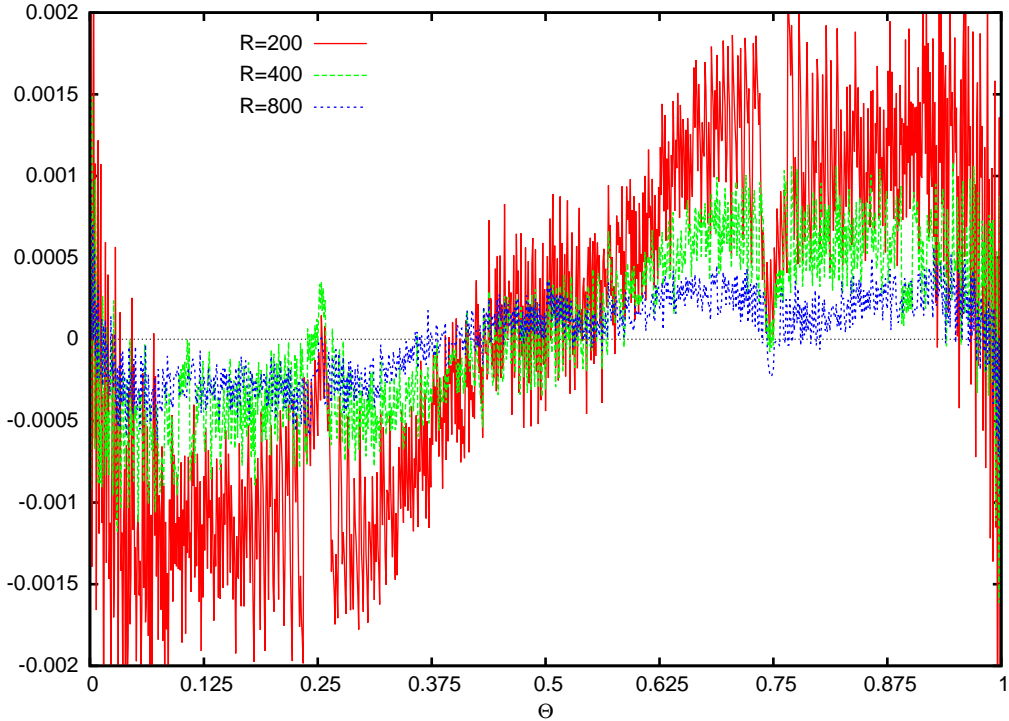}
\caption{The difference of the CDF's for the hitting 
distribution for the random 
walk and SLE$_6$ for the slit plane. The three curves are for 
$R=200, 400, 800$ to study the limit of the lattice spacing going to zero, 
i.e., $R \ra \infty$. }
\label{bad2_slit_exit_dif}
\end{figure}

In figure \ref{bad2_slit_exit_dif} we plot the difference between the 
simulation CDF and the SLE$_6$ CDF for the hitting distribution for the 
slit plane for radii of $200, 400,  800$. 
One sees that the difference is getting smaller as $R$ gets larger. 
The sawtooth nature of these plots can be understood as follows. 
In the simulation the random variable being studied is discrete. So its CDF
is a step function. When we subtract the smooth CDF from SLE$_6$ we get a 
sawtooth effect. As the radius gets larger this lattice effect gets 
smaller.

\begin{figure}[tbh]
\includegraphics{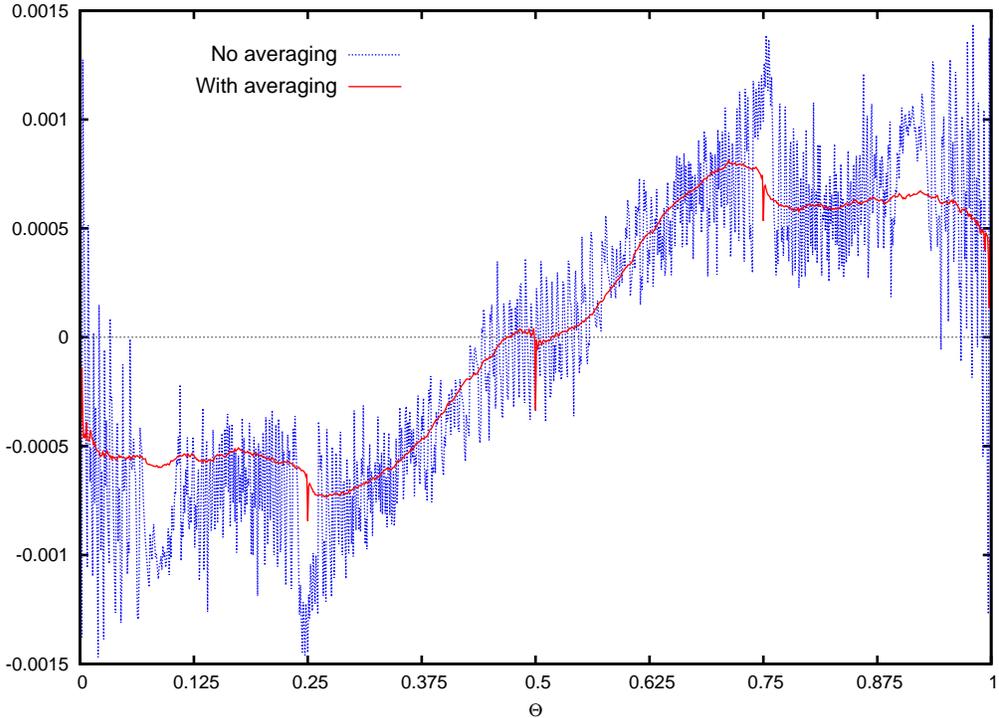}
\caption{The effect of averaging on the difference between the hitting
distributions for the random walk and SLE$_6$ in the slit plane. 
The highly oscillatory (blue) curve is for $R=375$. The smoother 
solid (red) curve averages $R$ over $[250,500]$. }
\label{bad2_smear_no_smear_exit}
\end{figure}

As discussed in the previous section, we can reduce 
this lattice effect by averaging over $R$.  
Figure \ref{bad2_smear_no_smear_exit} shows the effect of this averaging. 
There are two plots of the difference between the simulation CDF 
and the SLE$_6$ CDF. The sawtooth curve is for $R=375$ without averaging. 
For the other curve $R$ is averaged over $[250,500]$.
This plot shows that the averaging greatly reduces the sawtooth effect.
Note that the overall magnitude of the difference is not reduced by this 
averaging. The averaging smooths out the discreteness of the random variable,
but we still need to let the endpoints of the interval over which we average
go to infinity to take the scaling limit. 

\begin{figure}[tbh]
\includegraphics{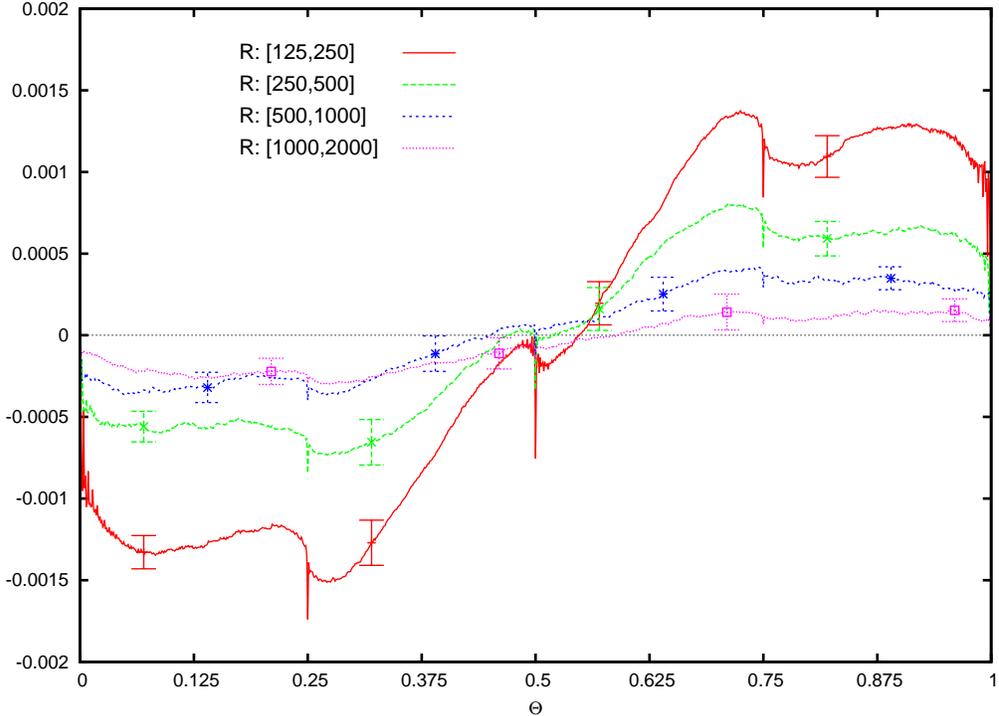}
\caption{The difference of the hitting distributions for the random walk and 
SLE$_6$ for the slit plane. $R$ is averaged over four different intervals
to study the limit of the lattice spacing going to zero, i.e., $R \ra \infty$. }
\label{bad2_slit_smear_exit_dif}
\end{figure}

\begin{figure}[tbh]
\includegraphics{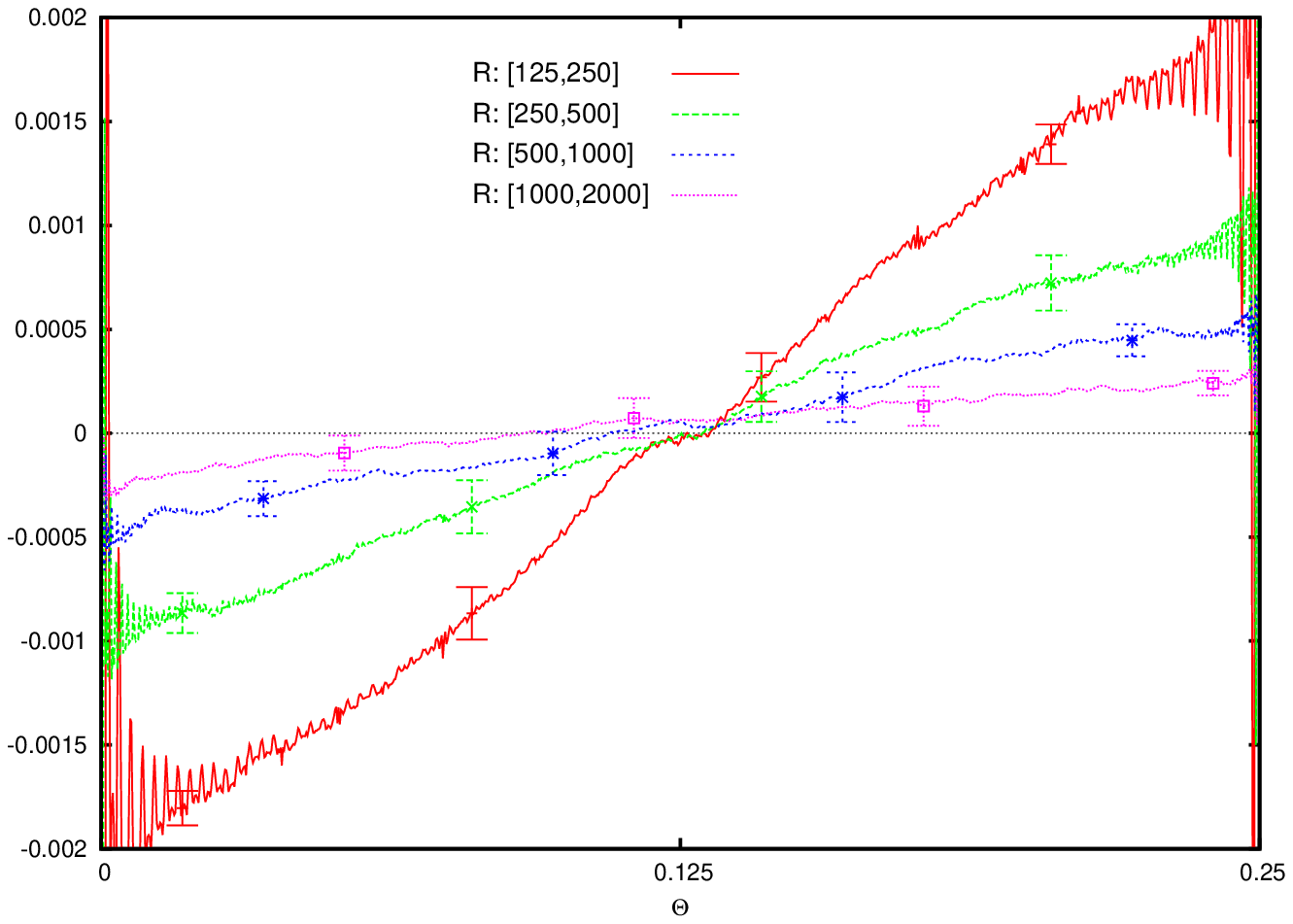}
\caption{The difference of the hitting distributions for the random walk and 
SLE$_6$ for the 90 degree wedge. $R$ is averaged over four different intervals
to study the limit of the lattice spacing going to zero, i.e., $R \ra \infty$. }
\label{bad2_90wedge_smear_exit_dif}
\end{figure}

\begin{figure}[tbh]
\includegraphics{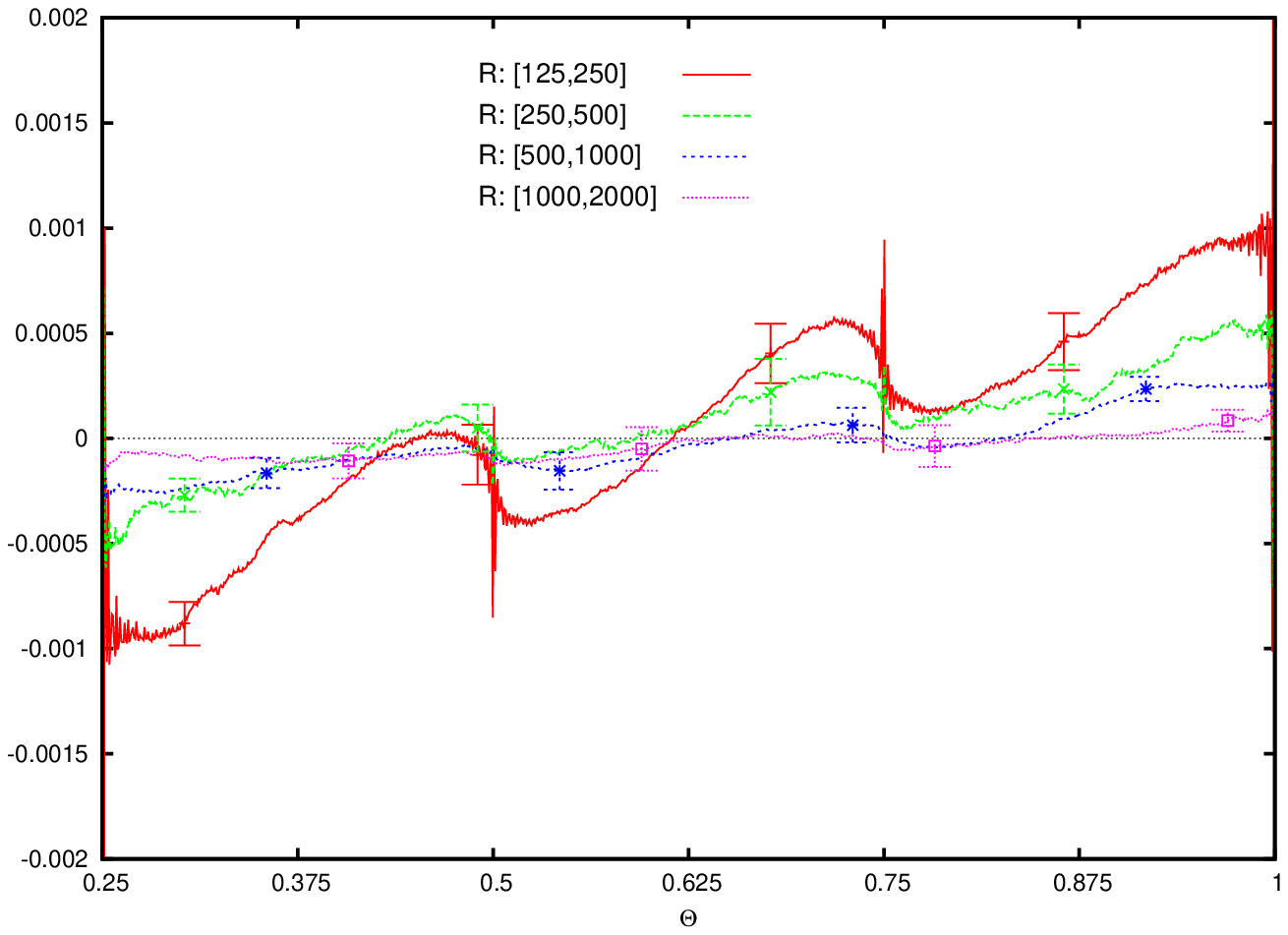}
\caption{The difference of the hitting distributions for the random walk and 
SLE$_6$ for the 270 degree wedge. 
$R$ is averaged over four different intervals to study 
the limit of the lattice spacing going to zero, i.e., $R \ra \infty$. }
\label{bad2_270wedge_smear_exit_dif}
\end{figure}

In figure \ref{bad2_slit_smear_exit_dif} we plot the difference between the 
simulation CDF and the SLE$_6$ CDF for the hitting distribution for the 
slit plane with this averaging. 
The intervals over which the radius $R$ is averaged are 
$[125,250], [250,500], [500,1000], [1000,2000]$.
In figures \ref{bad2_90wedge_smear_exit_dif} and 
\ref{bad2_270wedge_smear_exit_dif} we show the analogous plots for the 
90 degree wedge and the 270 degree wedge. 
In these plots the error bars shown represent plus/minus two standard 
deviations for the statistical errors that are a result of the number of 
samples being finite. The error bars do not represent the error that comes
from not having completely converged to the scaling limit. 
In all three figures the difference can be seen to decrease as $R$ increases
and is quite small. (Note the scale on the vertical axis.) 
For the largest value of $R$ the size of the difference is comparable to the 
statistical errors. 

\subsection{Pass right function}
 
We now study the probability of passing right of a point $R e^{i\theta}$ as a 
function of the polar angle $\theta$. 
As discussed before, $\delta=1/R$ is the effective lattice spacing and 
$n=N/R^{1/\nu}$ is the effective length of the walk (compared to the number 
of steps needed to reach the circle for the first time). 
To obtain the scaling limit we must let  $\delta$ go to zero 
and $n$ to $\infty$.  Since we average $R$ over an interval, both $\delta$ and 
$n$ vary over an interval in each simulation. When we give explicit values
for $\delta$ and $n$, we will compute them using the midpoint of the interval
for $R$. 

If we generate samples of walks with $N$ steps, then we can use these samples 
as samples of walks with $f N$ steps for $f<1$. (We do this for the sake of 
efficiency.) So in our simulations we have three parameters: $N_0, f$ and $r$.
$N_0$ is the number of steps in the walks being generated, while $N=fN_0$ is 
the number of steps that we actually use. The radius $R$ is taken to be
$R=rN_0^\nu$. So the effective lattice spacing $\delta$ and the length of 
the walk $n$ are given in terms of $N_0, f$ and $r$ by 
\bea
\delta=\frac{1}{r N_0^\nu}, \quad n=\frac{f N_0}{r^{1/\nu} N_0} 
= \frac{f}{r^{1/\nu}}
\eea

\begin{figure}[tbh]
\includegraphics{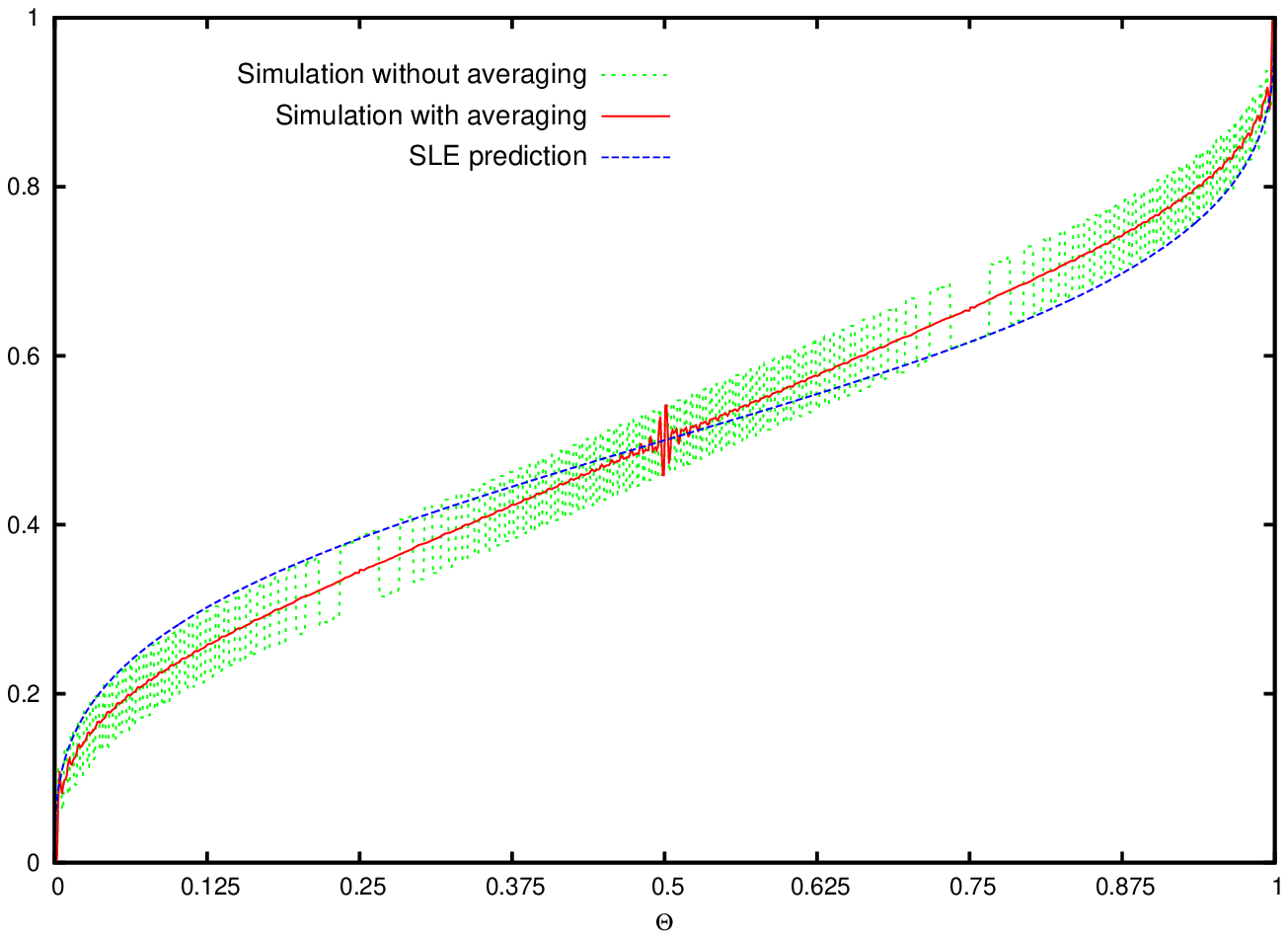}
\caption{The pass right function for the random walk and SLE$_6$. The highly
oscillatory (green) curve is the random walk without averaging. 
The solid (red) curve which passes though the center 
of the oscillatory curve is for the random walk with averaging. 
The dashed (blue) curve is for SLE$_6$. }
\label{bad2_smear_no_smear_pass_pfunc}
\end{figure}

To study the probability of passing right we use the averaging trick described
before. In figure \ref{bad2_smear_no_smear_pass_pfunc}
we plot the pass right function with $r$ averaged over $[0.05,0.1]$ 
and this function with $r=0.075$. For both of these plots $N=125K$. 
We also plot the SLE$_6$ prediction. 
Without averaging the CDF has large oscillations. The averaging removes these 
oscillations and the averaged CDF is centered in the middle of the 
oscillations. However, the averaged CDF does not agree at all with the 
SLE$_6$ CDF. 
The difference between the two curves is as large as $5\%$ for some $\theta$. 

\begin{figure}[tbh]
\includegraphics{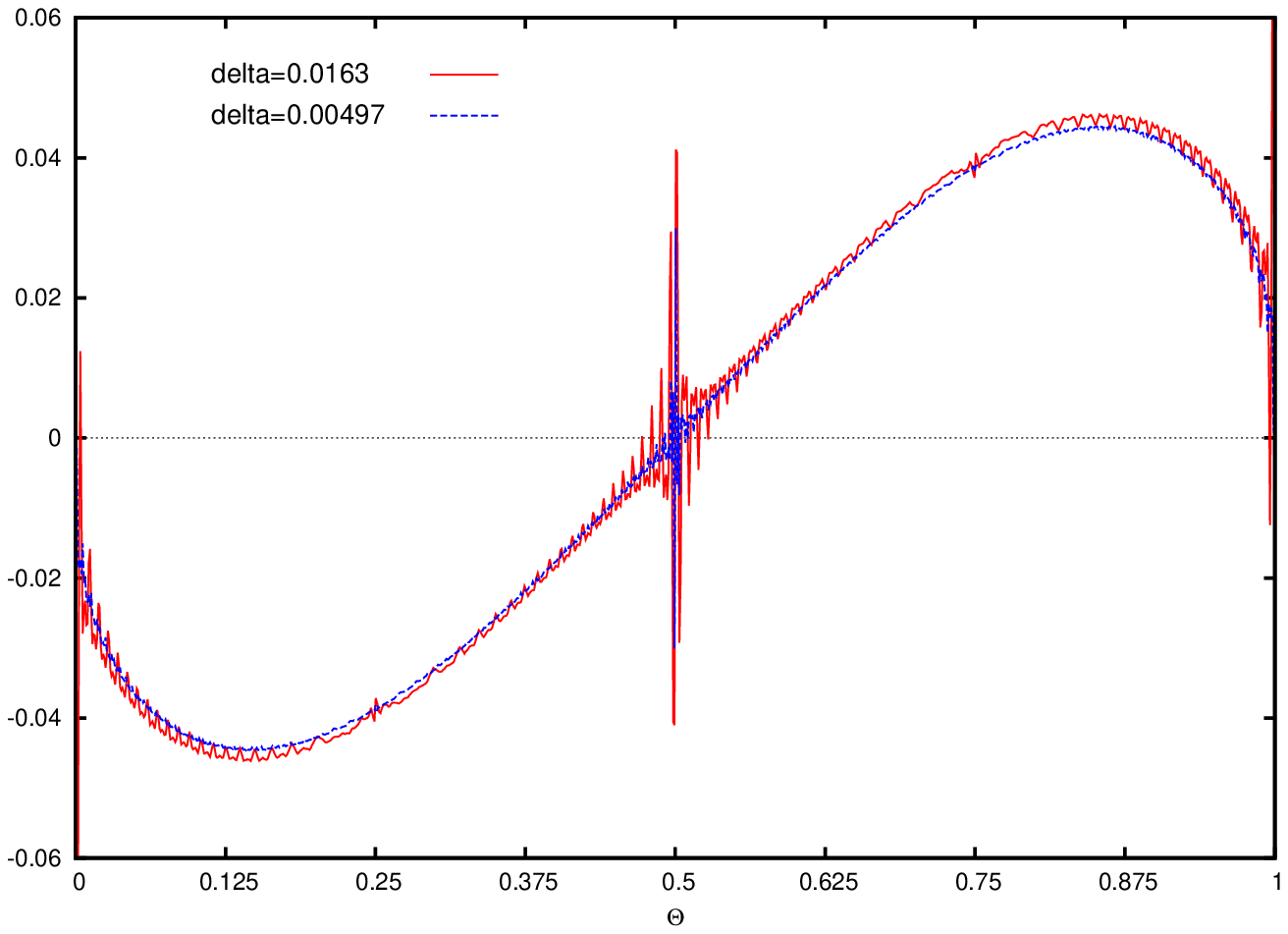}
\caption{The difference between the random walk and SLE$_6$ 
pass right functions for the slit plane when the length $n$ of the walk is 
kept fixed and the lattice spacing $\delta$ is varied.}
\label{bad2_slit_smear_pass_0.1_lattice}
\end{figure}

\begin{figure}[tbh]
\includegraphics{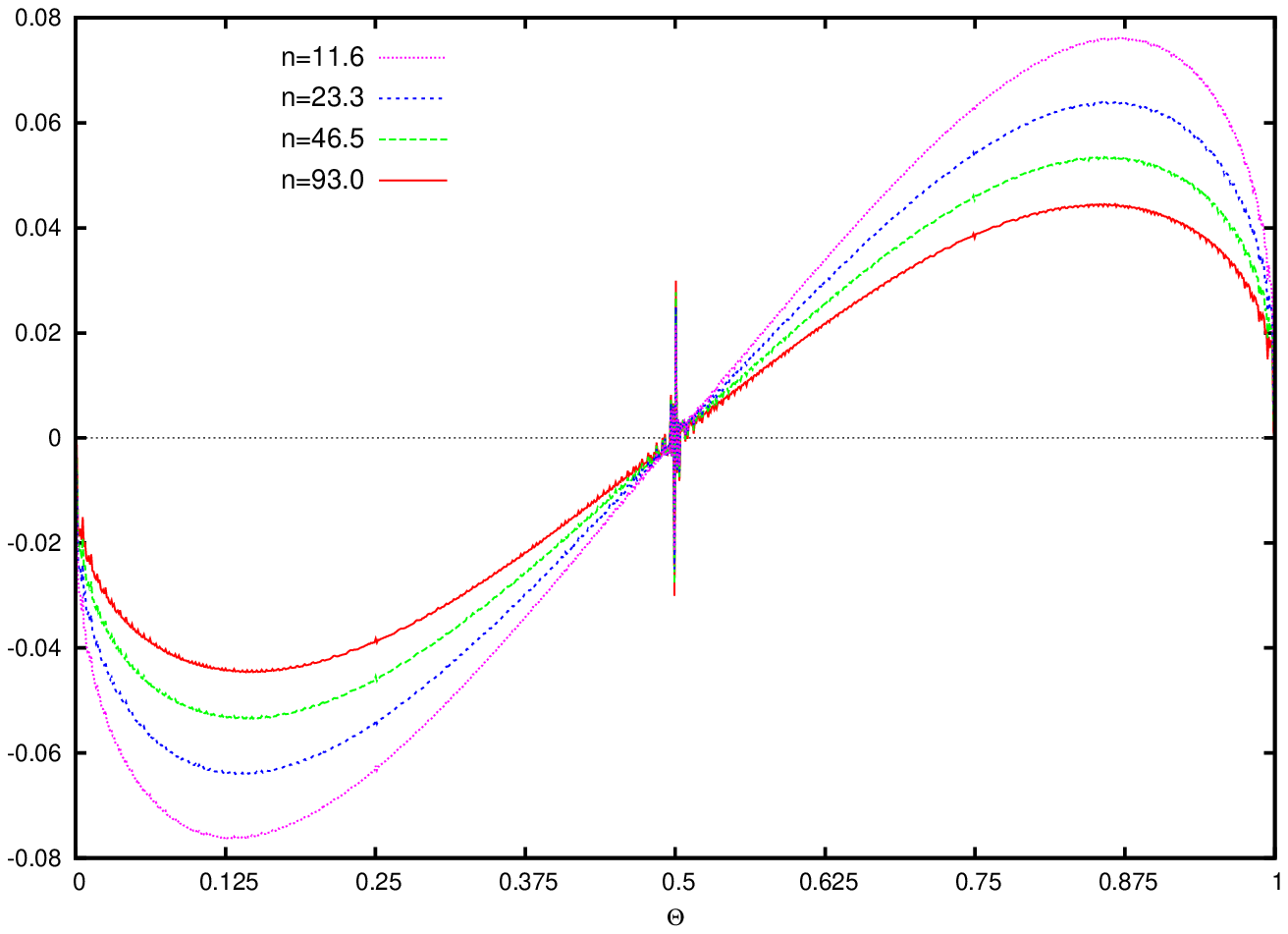}
\caption{The difference between the random walk and SLE$_6$ 
pass right functions for the slit plane when 
the lattice spacing $\delta$ is kept fixed and 
the length $n$ of the walk is varied.}
\label{bad2_slit_smear_pass_0.1}
\end{figure}

\begin{figure}[tbh]
\includegraphics{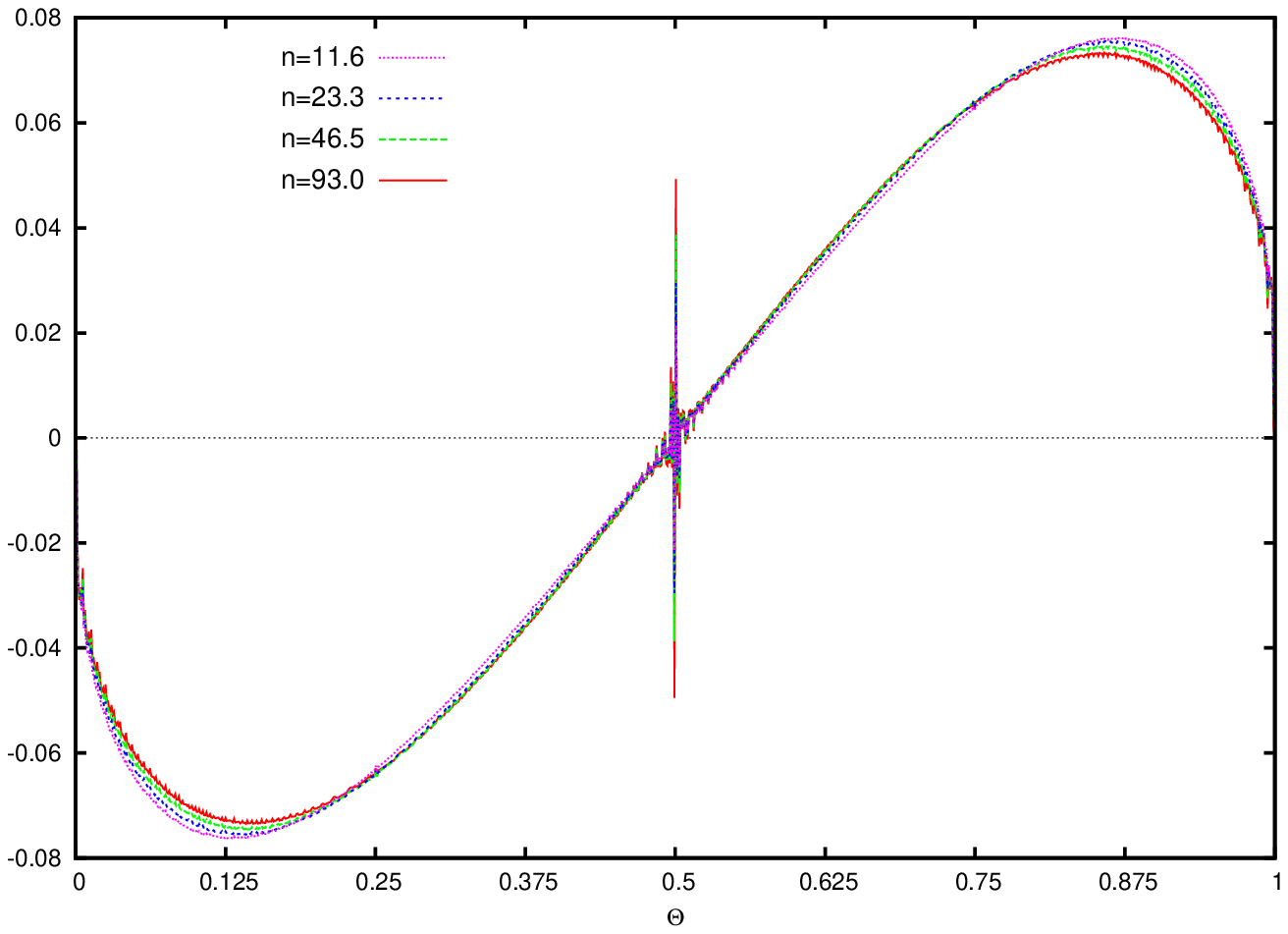}
\caption{The same differences shown in the previous plot are replotted 
but multiplied by a factor proportional to $n^p$ with $p=0.24$ to test if the 
differences in the previous plot are proportional to $n^{-p}$. }
\label{bad2_slit_smear_pass_rescaled_0.1_0.24}
\end{figure}

\begin{figure}[tbh]
\includegraphics{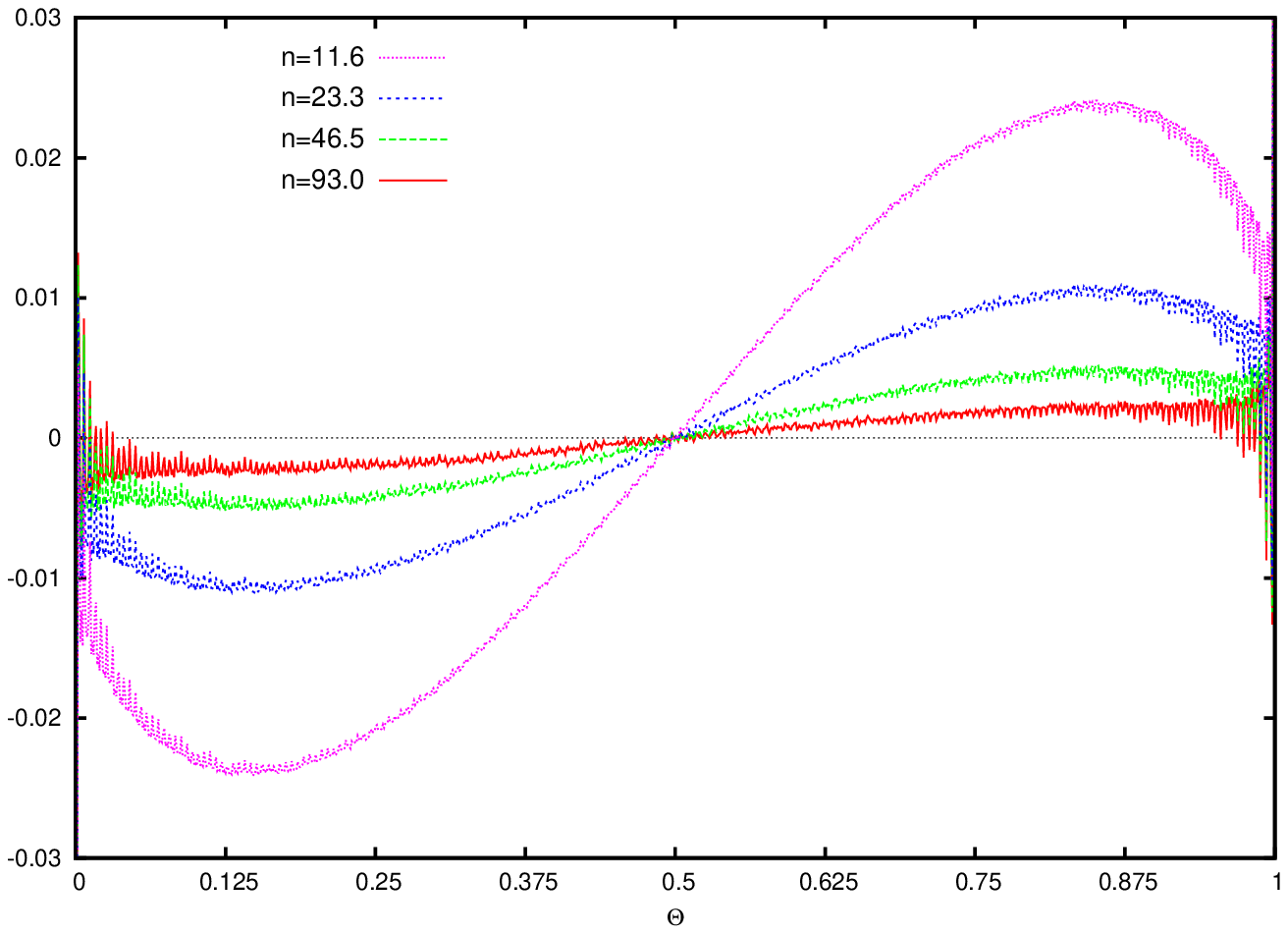}
\caption{The difference between the random walk and SLE$_6$ 
pass right functions for the 90 degree wedge when 
the lattice spacing $\delta$ is kept fixed and 
the length $n$ of the walk is varied.}
\label{bad2_90wedge_smear_pass_0.1}
\end{figure}

\begin{figure}[tbh]
\includegraphics{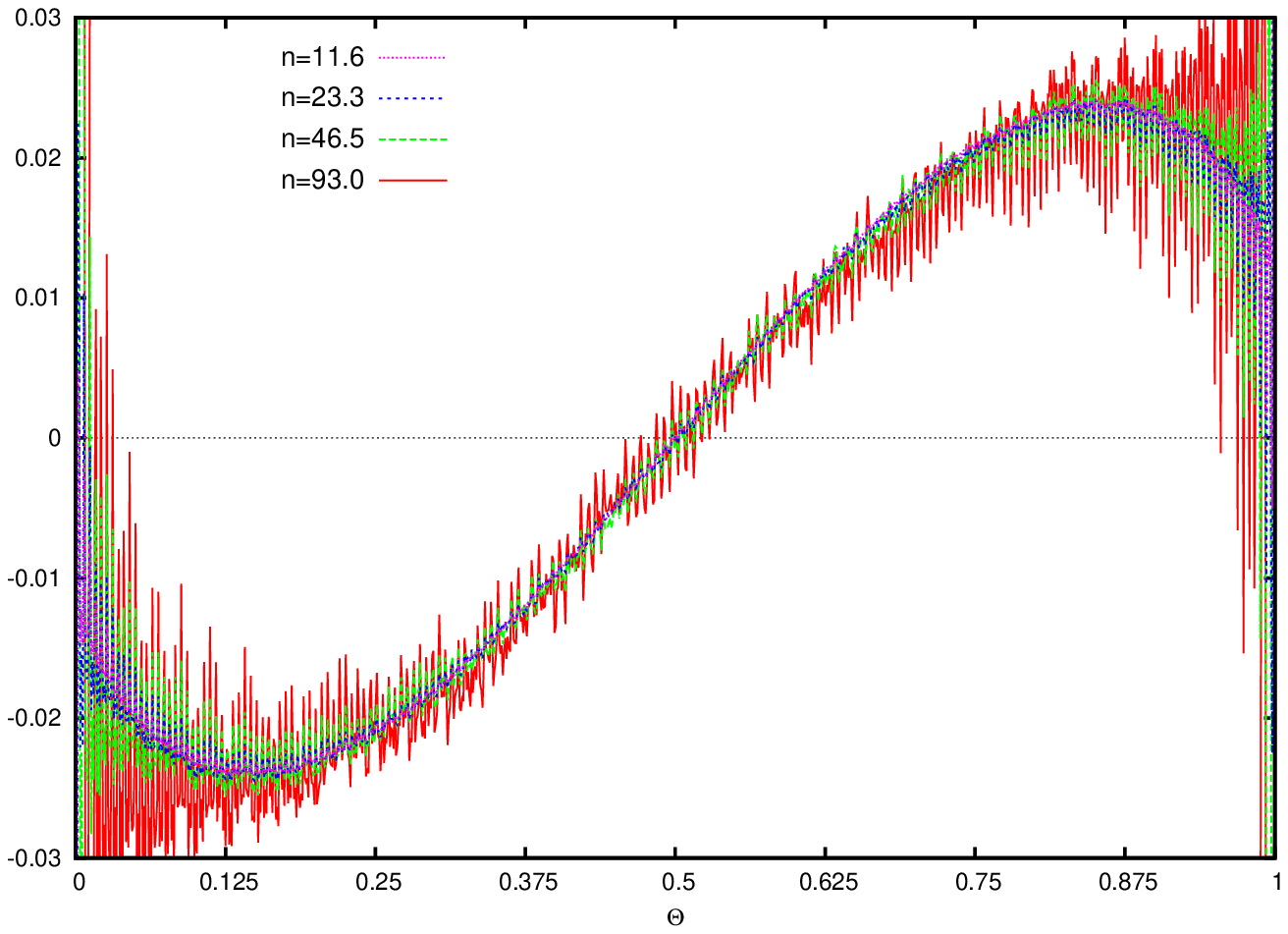}
\caption{The same differences shown in the previous plot are replotted 
but multiplied by a factor proportional to $n^p$ with $p=1.15$ to test if the 
differences in the previous plot are proportional to $n^{-p}$. }
\label{bad2_90wedge_smear_pass_rescaled_0.1_1.15}
\end{figure}

\begin{figure}[tbh]
\includegraphics{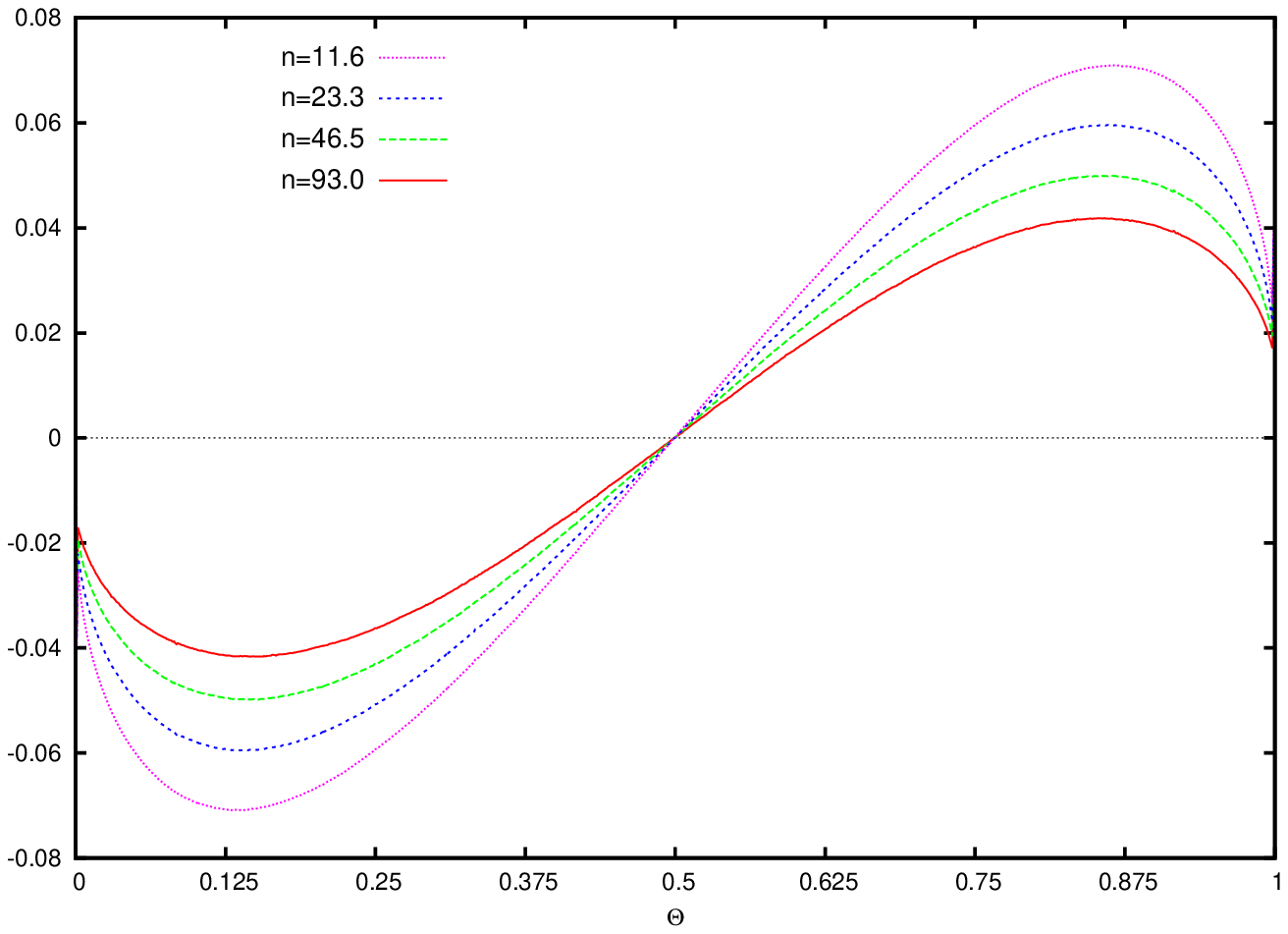}
\caption{The difference between the percolation explorer and SLE$_6$ 
pass right functions for the slit plane when 
the lattice spacing $\delta$ is kept fixed and 
the length $n$ of the percolation explorer is varied.}
\label{perc_slit_smear_pass_0.1}
\end{figure}

\begin{figure}[tbh]
\includegraphics{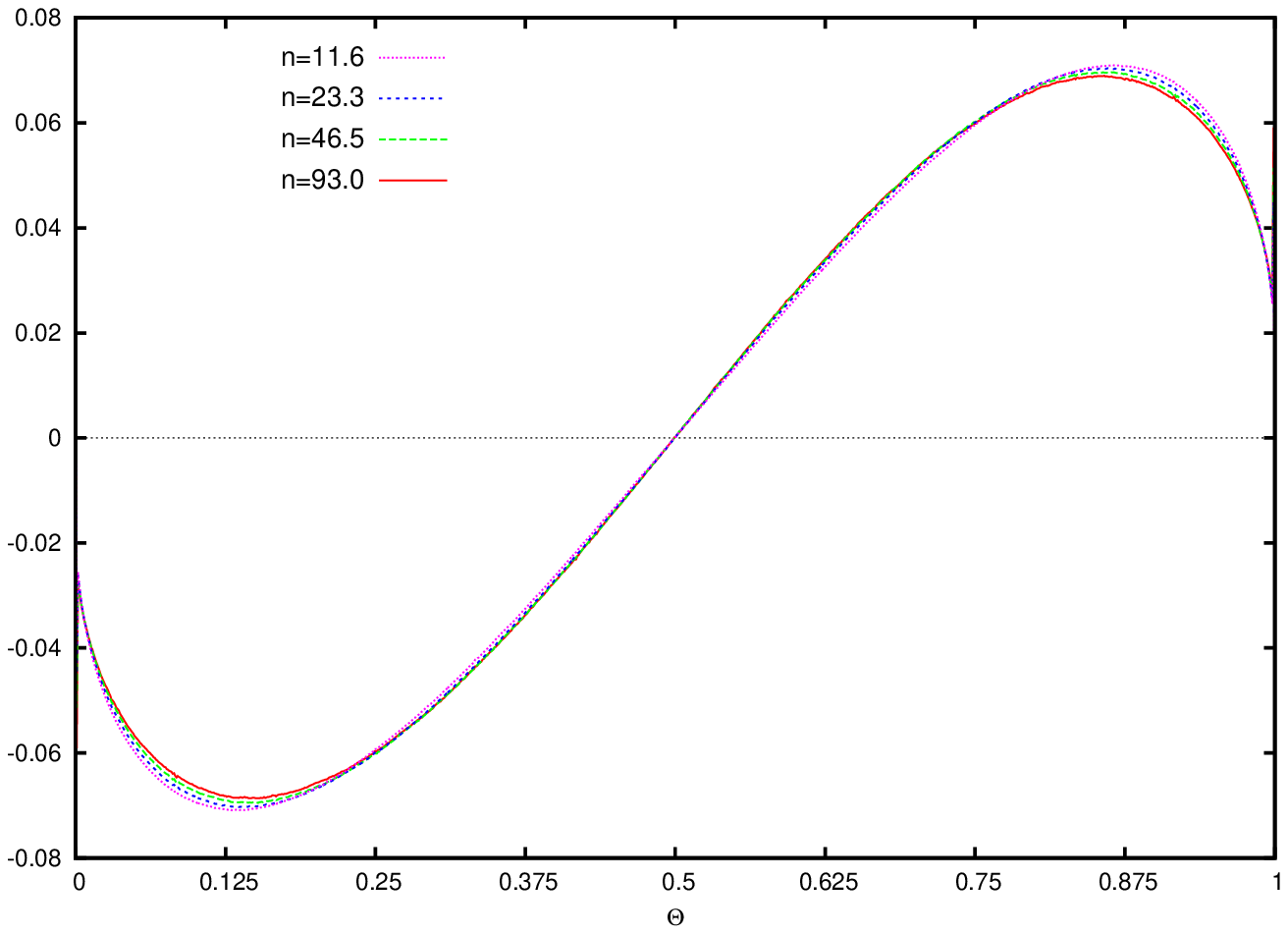}
\caption{The same differences for the percolation explorer 
shown in the previous plot are replotted but multiplied 
by a factor proportional to $n^p$ with $p=0.24$. }
\label{perc_slit_smear_pass_rescaled_0.1_0.24}
\end{figure}

The large discrepancy between the simulation and the SLE$_6$ prediction 
could be the result of the lattice spacing $\delta$ not being sufficiently 
small or the result of the length of the walk $n$ not being sufficiently 
large (or both). As we will show, it is primarily the result of $n$ not being 
sufficiently large. 
We first show that the effect of the lattice spacing $\delta$ on this large
discrepancy is negligible. Figure \ref{bad2_slit_smear_pass_0.1_lattice} 
shows two differences. Both differences use $r=0.1$ and $f=1$. 
One difference uses $N_0=125K$ and the other $N_0=1000K$. So the 
length of the walk $n$ is the same for the two differences. The effective 
lattice spacings are given by $\delta=1/(r N_0)^\nu$ which takes on the values
$0.0163$ and $0.00497$. So $\delta$ varies by a factor of 
$8^{4/7} \approx 3.3$, but the figure shows that the difference hardly changes. 
In this figure and the remaining figures in the paper we do not show any
error bars since they are tiny compared to the differences plotted in these 
figures. Two standard deviations are on the order of $10^{-4}$ which is 
much smaller than the vertical scale of the plots in figures 
\ref{bad2_slit_smear_pass_0.1_lattice} to 
\ref{perc_slit_smear_pass_rescaled_0.1_0.24}.

Next we fix $\delta$ and vary $n$. 
We do this by fixing $N_0=10^6$ and fixing the interval $r$ over which we 
average to be $[0.05,0.10]$. Then we vary $f$ over $1/8, 1/4, 1/2, 1$.
So in all four cases the lattice spacing $\delta$ is averaged over the same
interval while the length of the walk $n$ is averaged over 
four different intervals. We use the midpoint of the interval for $R$ 
to compute explicit values for $n$. They are $n=11.6, 23.3, 46.5, 93.0$. 
The resulting four differences are shown in 
figure \ref{bad2_slit_smear_pass_0.1} for the slit plane. 
The difference looks like it is going to zero as $n$ goes to 
infinity, and by studying the size of 
these curves it appears it is converging as $n^{-p}$ where a crude estimate
of $p$ is $0.24$. 

We test this by rescaling the difference cuves with 
a factor proportional to $n^p$. The curve for $n=93.0$ is rescaled by 
a factor of $2^{3p}$, the curve for $n=46.5$ by a factor of $2^{2p}$ and the 
curve for $n=23.3$ by a factor of $2^p$. The curve for $n=11.6$ 
is not rescaled. These rescaled curves are shown in figure
\ref{bad2_slit_smear_pass_rescaled_0.1_0.24} for $p=0.24$. 
They agree well. Note that all four differences in figure  
\ref{bad2_slit_smear_pass_0.1} use the same non-zero lattice spacing. So we 
do not expect them to converge exactly to zero. So we expect 
the rescaled differences in figure 
\ref{bad2_slit_smear_pass_rescaled_0.1_0.24} to differ slightly. 

In figure \ref{bad2_90wedge_smear_pass_0.1} we plot the difference 
for the 90 degree wedge with $\delta$ fixed and $n$ varying just as we 
did for the slit plane. (We take $N_0=10^6$,  average $r$ over $[0.05,0.10]$,
and vary $f$ over $1/8, 1/4, 1/2, 1$.)
For the 90 degree wedge the difference converges to zero much faster 
as $n$ goes to infinity. Again, the converence appears to go as $n^{-p}$, 
but now a crude estimate of $p$ gives the much larger value of $1.15$. 
In figure \ref{bad2_90wedge_smear_pass_rescaled_0.1_1.15} we rescale these
differences by a factor proportional to $n^p$ just as we did for the 
slit plane. The rescaled curves agree well. The lattice effects appear larger
here compared to the slit plane because the rescaling factors of $2^p, 2^{2p}$ 
and $2^{3p}$ are considerably larger than for the slit plane. 

While the simulations provide strong evidence that in the scaling limit
the pass right function for the random walk is converging to the SLE$_6$
prediction, the slow convergence of the $n \rightarrow \infty$ limit, 
especially for the slit plane, is rather surprising. 
So it is instructive to look at the analogous simulation for 
the percolation explorer on the hexagonal lattice
which is proven to converge to SLE$_6$. 
We give a brief description of this model and refer the reader to 
\cite{werner2007lectures} for more detail. 
We work in the half-plane. The hexagons above the horizontal 
axis are randomly colored black or white with equal probability. 
The hexagons along the horizontal axis that are to the right of the origin
are colored white, and those to the left of the origin are colored black. 
We now consider the interfaces between black and white hexagons. These 
interfaces will form loops with one exception. The choice of boundary 
condition forces there to be an interface which starts at the origin and 
is of infinite length. This random curve is the percolation explorer which
has been proved to converge to SLE$_6$. 
We have simulated this percolation explorer and computed the probability of 
passing right of the points along a circle just as we did for our random 
walk model. We have used the same parameters (and hence the same $\delta$ and 
$n$ values) that we did for the random walk model. 
In figure \ref{perc_slit_smear_pass_0.1} 
we plot the difference between the simulation CDF and the SLE$_6$ CDF
for the slit plane. So this figure is the analog of figure 
\ref{bad2_slit_smear_pass_0.1}. The difference curves from the 
percolation simulation are remarkably similar to the difference curves 
from our random walk simulations. 

Just as in figure \ref{bad2_slit_smear_pass_0.1} we see that the difference
is large, but appears to be going to zero. We can check this more 
carefully by rescaling. 
Figure  \ref{perc_slit_smear_pass_rescaled_0.1_0.24}
shows the rescaled differences for percolation,  so this figure is the analog
of figure \ref{bad2_slit_smear_pass_rescaled_0.1_0.24}.
We use the same power $p=0.24$ that we used for the random walk model. 

\subsection{Details of the simulation}

The algorithm to generate the walk is very simple. At each step we must 
determine which of the two possible steps are not allowed, i.e., which of these 
two sites have been visited before or belong to the boundary.
We check if they have been visited before using a hash table. 
If both sites have not been visited and are not on the boundary, 
we randomly choose between them with equal probability. 

The time to generate a walk is proportional to the number of steps. 
Some care is needed to be sure the time it takes to evaluate the pass 
right function does not dominate the computation time. In particular, 
computing the intergral involved in our averaging for each random 
walk sample would be prohibitively expensive. 
Instead, we incorporate the evaluation
of this integral into the Monte Carlo simulation. For each sample of the 
random walk we randomly choose an $R$ uniformly from $[R_0,R_1]$ and just
compute the pass right function or the hitting point for that radius $R$. 

In all our simulations we generate $10^8$ samples of the random walk. The
time needed depends on the domain and the parameters which control the number
of steps in the walk. Furthermore, these computations are done on 
a large cluster with many users. So we only have a crude estimate of the 
CPU time used. The longest simulations take 
roughly on the order of $200$ CPU-days to generate $10^8$ samples.

\section{Conclusions}

We have studied a random walk model on the Manhattan lattice which is 
not allowed to visit a site more than once. This model is known to be 
related to interfaces for bond percolation on a square lattice, 
so it is natural to conjecture that its scaling limit should be SLE$_6$. 
We have tested
this conjecture with Monte Carlo simulations of two quantities. 
One is a sort of hitting distribution. The other is the probability of 
passing to the right of a given point. 

For the hitting distribution we found excellent agreement with the 
SLE$_6$ prediction for this quantity. 
For the domains we consider, the scaling limit of the probability of passing 
right requires taking a double limit - letting the lattice spacing go to 
zero and the length of the walk to infinity. The convergence for the latter
is surprisingly slow and there are significant differences between the function
we find in our simulations and the SLE$_6$ prediction for this function. 
However, we can fit the difference quite well by a function proportional
to $n^{-p}$ where $n$ is the length of the walk. Moreover, the differences
we see in our model are remarkably similar to the analogous differences 
for the percolation explorer on the hexagonal lattice which been proven 
to converge to SLE$_6$. 
We conclude that our simulations provide strong support for the conjecture 
that the scaling limit of the random walk on the Manhattan lattice with no
self-intersections is SLE$_6$. 

\bigskip

{\it Acknowledgments:} 
This research was supported in part by NSF grant DMS-1500850.
An allocation of computer time from the UA Research Computing High Performance 
Computing (HPC) and High Throughput Computing (HTC) 
at the University of Arizona is gratefully acknowledged.

\end{document}